\documentclass[12pt,a4paper]{article}

\usepackage{latexsym}
\usepackage{amsmath,amssymb}
\usepackage{amsthm}
\usepackage[margin=3cm]{geometry}

\newtheorem*{teiri}{Theorem}

\newtheorem*{kei}{Corollary}

\title{On the duality and the derivation relations for multiple zeta values}
\author{Naho Kawasaki and Tatsushi Tanaka}
\date{}

\begin{document}

\maketitle
\thispagestyle{plain}
 
\begin{abstract}
We consider the problem of deducing the duality relation from the extended double shuffle relation for multiple zeta values. Especially we prove that the duality relation for double zeta values and that for the sum of multiple zeta values whose first components are 2's are deduced from the derivation relation, which is known as a subclass of the extended double shuffle relation. 
\end{abstract}


\section{Introduction}

The multiple zeta values (MZVs for short) are defined for positive integers $k_1,k_2,\ldots,k_n$ with $k_1 \geq 2$ by the convergent series
\begin{eqnarray*}
\zeta(k_1,k_2,\ldots,k_n)
 =  \sum_{m_1>m_2> \cdots >m_n>0}\frac{1}{{m_1}^{k_1}{m_2}^{k_2} \cdots {m_n}^{k_n}} . 
\end{eqnarray*}
We call $k=k_1+k_2+\cdots+k_n$ weight, $n$ depth, and $\#\{i|k_i>1 \ (1\leq i \leq n)\}$ height. It is important to state many identities among MZVs explicitly to clarify the algebraic structures among them. 

It is conjectured that the extended double shuffle relations \cite{ikz} give all linear  relations among MZVs. In \cite{ikz}, K.~Ihara, M.~Kaneko and D.~Zagier showed that the sum formula(\cite{gr,z2}), Hoffman's relation(\cite{hof1}), and the derivation relation(\cite{ikz}) for MZVs are subject to the extended double shuffle relation. Also, they showed   the extended double shuffle relation together with the duality relation implies Ohno's relation(\cite{o}). However, it is still not known whether the duality relation itself is implied by the extended double shuffle relation or not. 

In \cite{kaj}, J. Kajikawa showed the duality relation for the sum of MZVs with fixed weight, depth, and height is deduced from the derivation relation (and hence from the extended double shuffle relation). In this paper, we observe another lineage, that is the duality for the sum with fixed weight, depth, and $k_1$. In particular we here show that the duality relations of the following two cases are also implied by the derivation relation:
\begin{itemize}
\item the case for the double zeta values, 
\item the case for the sum of MZVs with fixed weight, depth, and $k_1=2$. 
\end{itemize}
 

\section{Algebraic setup and main theorem}
\subsection{Algebraic setup}

To state our main theorem, we use the algebraic setup introduced by Hoffman \cite{hof2}, which encodes the MZVs as monomials in the noncommutative polynomial ring $\mathfrak{h}=\mathbb{Q} \langle x,y \rangle$. (Also see \cite{ikz,ta2} for example.) We denote by $\{x,y\}^*$ the set of words over the alphabet $\{x,y\}$. Let ${\mathfrak{h}}^0 = \mathbb{Q}+x \mathfrak{h} y$, a subring of $\mathfrak{h}$. Define the $\mathbb{Q}$-linear map $Z : {\mathfrak{h}}^0 \rightarrow \mathbb{R}$ by 
\begin{eqnarray*}
Z(1) & = & 1, \\
Z(x^{k_1 -1}y x^{k_2 -1}y \cdots x^{k_n -1} y) & = & \zeta(k_1,k_2, \ldots ,k_n) \quad (k_1>1). 
\end{eqnarray*}
The weight $k$ and the depth $n$ of $\zeta(k_1,k_2, \ldots ,k_n)$ correspond to the total degree and the degree in $y$ of the monomial $x^{k_1 -1}y x^{k_2 -1}y \cdots x^{k_n -1} y$, respectively. To obtain a linear relation for MZVs is nothing but to obtain an element in the kernel of $Z$. \\

$\bullet$ (Duality relation). 
Let $\tau$ be the anti-automorphism on $\mathfrak{h}$ defined by
\begin{eqnarray*}
\tau(x)=y, \quad \tau(y)=x. 
\end{eqnarray*}
The map $\tau$ preserves ${\mathfrak{h}}^0$. The duality relation is stated as
\begin{equation*}
(1-\tau)(w) \in \ker Z
\end{equation*}
for any $w \in \mathfrak{h}^0$. \\

$\bullet$ (Derivation relation). 
For a positive integer $n$, we define the $\mathbb{Q}$-linear map ${\partial}_n:\mathfrak{h} \rightarrow \mathfrak{h}$ by 
\begin{equation*}
{\partial}_n (x) = -{\partial}_n (y) = x (x+y)^{n-1}y 
\end{equation*}
and Leibniz rule
\begin{equation*}
{\partial}_n (v w) = {\partial}_n (v) w + v {\partial}_n (w) \ \ \ \ (v, w \in \mathfrak{h}).
\end{equation*}
We find in \cite{ikz} that $\mbox{Im} {\partial}_n \subset \mathfrak{h}^0$ and $\deg({\partial}_n(w))=\deg w+n$ for $w\in\{x,y\}^*$. 
We also find that $\partial_n$ and $\partial_m$ commute for any $n,m\geq 1$. 
The derivation relation is then stated as  
\begin{equation*}
{\partial}_n (w) \in \ker Z
\end{equation*}
for any positive integer $n$ and any $w \in {\mathfrak{h}}^0$. \\

Let $\Theta$ be the map on $\hat{\mathfrak{h}}= \mathbb{Q} \langle \langle x,y \rangle\rangle$, the completion of $\mathfrak{h}$, defined by 
\begin{eqnarray*}
\Theta = {\rm{exp}} \left( \sum_{n \geq 1} \frac{{\partial}_n}{n} \right). 
\end{eqnarray*}
This $\Theta$ is an element in the graded ring $\mathbb{Q} [[{\partial}_1, {\partial}_2,\ldots]]$, where the degree of the map ${\partial}_n$ is given as $n$. 
Let ${\theta}_l$ be the degree $l$ homogeneous part of the map $\Theta$. 
For example, we have
\begin{eqnarray*}
{\theta}_1= {\partial}_1, \ \ 
{\theta}_2= \frac{1}{2} ({\partial}_2 + {{\partial}_1}^{2}), \ \ 
{\theta}_3= \frac{1}{6} (2 {\partial}_3 + 3 {\partial}_2 {\partial}_1 + {{\partial}_1}^{3}), \ldots. 
\end{eqnarray*}
We extend $\theta_l$ naturally to the map on $\hat{{\mathfrak{h}}}$. 


\subsection{Main theorem}
Our main theorem is then stated as follows. 
\begin{teiri}\label{teiri1}
The following two identities hold in $\hat{{\mathfrak{h}}^0}$, the completion of ${\mathfrak{h}}^0$. 
\begin{description}
\item[(i)]
For any positive integer $m$, we have 
\begin{eqnarray*}
&  & (1- \tau) \left( x^m y \frac{1}{1-x} y  \right) \\
& = & (\Theta-1) \left( x^m y \left( 1-\frac{1}{1-x} y \right) \right) \\
&   & - \sum_{i=1}^{m-1}{\theta}_i \left (x^{m-i} y
+ \left( x \frac{1}{1-y} \right) ^{m-i} xy 
- \left( x \frac{1}{1-y} \right) ^{m-i-1} xy \right) .
\end{eqnarray*}
\item[(ii)]
For any positive integer $n$, we have 
\begin{eqnarray*}
&   & (1- \tau) \left( xy \left( \frac{1}{1-x} y \right) ^{n-1} \right) \\
& = & (\Theta -1) \left( x \frac{1- x^{n-1}}{1-x} y \left( 1-\frac{1}{1-x} y \right) \right) - \sum_{l=1}^{n-2} {\theta}_l \left( x \frac{1- x^{n-l-1}}{1-x} y \left( 1- \frac{1}{1-x} y \right) \right).  
\end{eqnarray*}
\end{description}
Here, the empty sum is regarded as 0.
\end{teiri}

These formulas are more explicit than the following corollary as asserted in $\S1$. 
\begin{kei}
\begin{description}
\item[(i)] For any positive integer $s$ and non-negative integer $t$, we have
\begin{equation*}
(1-\tau)(x^s y x^t y) \in \sum_{n \geq 1} \partial_n (\mathfrak{h}^0).
\end{equation*}
\item[(ii)] For any positive integers $s,t$ with $s > t \geq 1$, we have
\begin{equation*}
(1-\tau) \left( \sum_{\substack{w \in \mathfrak{h}^1, \ \deg w=s-2, \\ \deg_y(w)=t-1}} xyw \right) \in \sum_{n \geq 1} \partial_n (\mathfrak{h}^0).
\end{equation*}
Here $\deg_y(w)$ denotes the degree of $w$ in $y$.  
\end{description}
\end{kei}

Corollary (i) shows that the duality relation for each double zeta value can be viewed as a sort of the derivation relation. 
Corollary (ii) shows that the duality relation for the sum of MZVs with fixed weight, depth, and $k_1=2$ can be viewed as a sort of the derivation relation. 
We remark that, for example, the dualities for $\zeta(3,2)$ and $\zeta(2,3)$ (i.e. $\zeta(3,2)=\zeta(2,2,1)$ and $\zeta(2,3)=\zeta(2,1,2)$) are regarded as derivation relations because of our theorem (ii), though we cannot conclude them by Kajikawa's result \cite{kaj}. 

Each row, except for the first one which gives division in weights, of the following table shows the numbers of linearly independent relations for the class explained shortly (and assigned numbering from 1. to 7.) in the leftmost box. For example, no.1 gives the numbers of linearly independent relations of the duality for the sum of fixed weight, depth, and height, no.3 gives those of the union class of no.1 and no.2, and no.7 gives those of the intersection of no.4 and no.5. The class no.6 is known to be equivalent to Ohno relation \cite{o}. These numbers are computed by using Risa/Asir, an open source general computer algebra system. 
\begin{table}[h]\hspace{-15pt}
\begin{footnotesize}
\begin{tabular}{|l||c|c|c|c|c|c|c|c|c|c|c|c|} \hline
wt & 3 & 4 & 5 & 6 & 7 & 8 & 9 & 10 & 11 & 12 & 13 & 14  \\ \hline \hline
1. Duality (fixed wt, dep, and ht) & 1 & 1 & 3 & 3 & 6 & 6 & 10 & 10 & 15 & 15 & 21 & 21 \\ \hline
2. Duality (fixed wt, dep, and $k_1$) & 1 & 1 & 4 & 6 & 11 & 15 & 22 & 28 & 37 & 45 & 56 & 66  \\ \hline
3. 1. $\cup$ 2.  & 1 & 1 & 4 & 6 & 12 & 16 & 25 & 31 & 43 & 51 & 66 & 76  \\ \hline
4. Duality & 1 & 1 & 4 & 6 & 16 & 28 & 64 & 120 & 256 & 496 & 1024 & 2016 \\ \hline
5. Derivation relation & 1 & 2 & 5 & 10 & 22 & 44 & 90 & 181 & 363 & 727 & 1456 &  \\ \hline
6. 4. $\cup$ 5. (Ohno relation) & 1 & 2 & 5 & 10 & 23 & 46 & 98 & 199 & 411 & 830 & 1691 & \\ \hline
7. 4. $\cap$ 5. & 1 & 1 & 4 & 6 & 15 & 26 & 56 & 102 & 208 & 393 & 789 &  \\ \hline
\end{tabular}
\end{footnotesize}
\end{table}

\noindent We confirm that the space of the duality for the sum with fixed weight, depth, and $k_1$ (no.2) is included in the space of the derivation relation (no.5) if the weight is less than or equal to 12 due to an experiment using a computer, but in general it is still unknown. More precisely, we conjecture that
\begin{equation*}
 (1-\tau)\left(x^{m-1}y \left( \frac{1}{1-x}y \right) ^{n-1}\right)\stackrel{?}{\in}\sum_{n\geq 1}\partial_n(\mathfrak{H}^0)
\end{equation*}
for any $m,\,n\geq 3$. The sequence in the last line (no.7) is obtained just by the inclusion-exclusion property, that is no.4 plus no.5 minus no.6 for each column. We see that the sequence in no.3 is far smaller than the sequence in no.7. How can we explain the sequence in no.7 concretely?

\section{Proof of the theorem}
Let $\mathfrak{h}[[u]]$ be the formal power series ring over $\mathfrak{h}$ with an indeterminate $u$, and ${\Delta}_u$ the automorphism of $\mathfrak{h}[[u]]$ whose images of generators are 
\begin{equation*}
{\Delta}_u (u) = u, \ \ \ \ {\Delta}_u (x) = x \frac{1}{1-yu}, \ \ \ \ {\Delta}_u (y)= (1-xu-yu) \frac{y}{1-yu}.
\end{equation*}
The images of the generators of the inverse map ${{\Delta}_u}^{-1}$ are then given by 
\begin{equation*}\label{eq2}
{{\Delta}_u}^{-1} (u) = u, \ \ \ \ {{\Delta}_u}^{-1} (x) = \frac{x}{1-xu} (1-xu-yu), \ \ \ \ {{\Delta}_u}^{-1} (y)= \frac{1}{1-xu} y.
\end{equation*}
(Also see \cite{ikz, kaj, ta1}.)  In fact, we have $\Delta_u=\sum_{l=0}^{\infty} \theta_l u^l$. 
Putting $u=1$, we see that the map ${\Delta}_1$ (which is naturally extended to as the automorphism of $\hat{\mathfrak{h}}$, the completion of $\mathfrak{h}$) coincides with $\Theta$ on ${\hat{\mathfrak{h}}}^0$. Hence, we may extend $\Theta$ to the automorphism of $\hat{\mathfrak{h}}$ via the identification with $\Delta_1$. 
\begin{proof}
\begin{description}
\item[(i)] We  calculate the generating function of both sides of identity (i). 
We multiply both sides by $u^m$ and add them up for $m \geq 1$. 
Calculating generating function of LHS of (i), we obtain
\begin{equation}\label{eq12}
\sum_{m \geq 1} \left( x^m y \frac{1}{1-x} y -x \frac{1}{1-y}x y^m \right) u^m 
 =  \frac{xu}{1-xu} y \frac{1}{1-x} y - x \frac{1}{1-y} x \frac{yu}{1-yu}. 
\end{equation}
On the other hand, calculating generating function of the first term of RHS of (i), we obtain
\begin{eqnarray}
&  & \sum_{m \geq 1} (\Theta-1) \left( x^m y \left( 1-\frac{1}{1-x} y \right) \right) u^m \nonumber\\
& = & (\Theta-1) \left( \frac{xu}{1-xu} y \left( 1-\frac{1}{1-x} y \right) \right) \nonumber\\
\label{eq8}
& = & xu \frac{1}{1-xu-y} (1-x-y)y - \frac{xu}{1-xu} y \left( 1- \frac{1}{1-x}y \right). 
\end{eqnarray}
And for the second term of RHS of (i), we obtain
\begin{eqnarray}
- \sum_{m \geq 1} \sum_{i=1}^{m-2} {\theta}_i (x^{m-i} y) u^m 
& = & - ({\Delta}_u -1) \left( \frac{(xu)^2}{1-xu} y \right) \nonumber\\
\label{eq9}
& = & -xu \frac{1}{1-yu} x \frac{yu}{1-yu} + \frac{xu}{1-xu} y -xyu. 
\end{eqnarray}
And for the third term of RHS of (i), we obtain
\begin{eqnarray}
&  & - \sum_{m \geq 1} \sum_{i=1}^{m-1} {\theta}_i \left( \left( x \frac{1}{1-y} \right) ^{m-i} xy \right) u^m \nonumber\\
\label{eq10}
& = & - {\Delta}_u \left( xu \frac{1}{1-xu-y} xy \right) + xu \frac{1}{1-xu-y} xy.
\end{eqnarray}
And for the fourth term of RHS of (i), we obtain
\begin{eqnarray}
&  & \sum_{m \geq 1} \sum_{i=1}^{m-2} {\theta}_i \left( \left( x \frac{1}{1-y} \right) ^{m-i-1} xy \right) u^m \nonumber\\
\label{eq11}
& = & u {\Delta}_u \left( xu \frac{1}{1-xu-y} xy \right) - xu \frac{1}{1-xu-y} xyu.
\end{eqnarray}
Adding the second term of (\ref{eq8}) and the second term of (\ref{eq9}), we obtain
\begin{eqnarray*}
 - \frac{xu}{1-xu} y \left( 1- \frac{1}{1-x}y \right) + \frac{xu}{1-xu} y 
 = \frac{xu}{1-xu} y \frac{1}{1-x}y. 
\end{eqnarray*}
This is the first term of (\ref{eq12}). 
Using 
\begin{equation*}
\frac{1}{1-xu-y}=\frac{1}{1-\frac{1}{1-xu}y} \cdot \frac{1}{1-xu},
\end{equation*}
one has 
\begin{equation*}
 x \frac{1}{1-y} x \frac{yu}{1-yu} - xu \frac{1}{1-yu} x \frac{yu}{1-yu} = (1-u) {\Delta}_u  \left( xu \frac{1}{1-xu-y} xy \right), 
\end{equation*}
 which means the sum of the first terms of (\ref{eq9}), (\ref{eq10}) and (\ref{eq11}) turns into the second term of (\ref{eq12}). 
The remaining terms in (\ref{eq8}), (\ref{eq9}), (\ref{eq10}), and (\ref{eq11}) cancel. 
This shows that the generating functions of the both sides are equal and therefore we conclude (i) in the theorem. 


\item[(ii)] Similary we calculate the generating functions of both sides of identity (ii).  We multiply both sides by $u^{n-1}$ and add them up for $n \geq 1$. Calculating generating function of LHS of (ii), we obtain
\begin{eqnarray}
&  & \sum_{n \geq 1} \left(  xy \left( \frac{1}{1-x} y \right) ^{n-1} - \left( x \frac{1}{1-y} \right) ^{n-1} xy \right) u^{n-1} \nonumber\\
\label{eq3}
& = & xy \frac{1}{1-x-yu} (1-x) - (1-y) \frac{1}{1-xu-y} xy. 
\end{eqnarray}
On the other hand, we calculate 
\begin{equation}\label{eq1}
\Theta \left( x \frac{1- x^{n-1}}{1-x} y \left( 1-\frac{1}{1-x} y \right) \right) - \sum_{l=0}^{n-2} {\theta}_l \left( x \frac{1- x^{n-l-1}}{1-x} y \left( 1- \frac{1}{1-x} y \right) \right)
\end{equation}
which is equal to generating function of RHS of (ii). 
Calculating generating function of the first term of (\ref{eq1}), we obtain
\begin{eqnarray}\label{eq4}
&    & \sum_{n \geq 1} \Theta \left( x \frac{1-x^{n-1}}{1-x} y \left( 1- \frac{1}{1-x} y \right) \right) u^{n-1} \nonumber\\
& = & xy\frac{1}{1-u} - x\frac{1}{1-xu-y}(1-y)y \nonumber\\
& = & xy \frac{1}{1-u} - (1-y) \frac{1}{1-y-xu} xy. 
\end{eqnarray}
And for the second term of  (\ref{eq1}), we obtain
\begin{eqnarray}
&  & -  \sum_{n \geq 1} \sum_{l=0}^{n-2} {\theta}_l \left( x \frac{1-x^{n-l-1}}{1-x} y \left( 1-\frac{1}{1-x}y \right) \right) u^{n-1} \nonumber\\
& = & - {\Delta}_u \left( \frac{xu}{1-xu} y \left( 1-\frac{1}{1-x}y \right) \right) \frac{1}{1-u}
 \nonumber\\
\label{eq6}
& = & - xy \frac{u}{1-yu} \frac{1}{1-u} + xy \frac{1}{1-x-yu} (1-xu-yu) \frac{yu}{1-yu} \frac{1}{1-u} \nonumber \\
& = & xy\frac{u}{1-u}+xyu\frac{1}{1-x-yu}y
\end{eqnarray}
In this case, we can show that (\ref{eq3}) is equal to the sum of (\ref{eq4}) and (\ref{eq6}) by straightforward calculation. Therefore we conclude (ii) in the theorem. 
\end{description}
\end{proof}

\section*{Acknowledgement}
The authors express great thanks to Prof. Yasuo Ohno and Prof. Masanobu Kaneko for helpful comments and advice in the MZV conference at Kinki university on Feb., 2017. This work is supported in part by Kyoto Sangyo University Research Grants (2016). 



{\sc Mathematical Institute, Tohoku University, Sendai 980-8578, Japan.}

{\it E-mail address}: {\tt naho.kawasaki.p7@dc.tohoku.ac.jp} \\

{\sc Department of Mathematics, Faculty of Science, Kyoto Sangyo University, Motoyama Kamigamo Kita-Ku Kyoto-City, 603-8555, Japan.}

{\it E-mail address}: {\tt t.tanaka@cc.kyoto-su.ac.jp}

\end{document}